\documentclass{amsart}

\usepackage[utf8]{inputenc}

\usepackage{boxedminipage}
\usepackage{amsfonts}
\usepackage{amsmath} 
\usepackage{amssymb}
\usepackage{graphicx}
\usepackage{amsthm}
\usepackage{t1enc}
\usepackage{subfig}
\usepackage{xfrac}

\usepackage{todonotes}

\newtheorem{theorem}{Theorem}

\newtheorem{example}[theorem]{Example}

\title[On free semigroups of affine maps on the real line]{On free semigroups of affine maps \\ on the real line}

\author{Alexander Kolpakov}
\address{\parbox{\linewidth}{Institut de Math\'ematiques, Universit\'e de Neuch\^atel, 2000 Neuch\^atel, Switzerland \\
Laboratory of combinatorial and geometric structures, Moscow Institute of Physics \\ and Technology, Dolgoprudny, Russia}}
\email[]{kolpakov (dot) alexander (at) gmail (dot) com}

\author{Alexey Talambutsa}
\address{\parbox{\linewidth}{Steklov Mathematical Institute of RAS, 8 Gubkina St., 119991 Moscow, Russia\\
HSE University, 11 Pokrovsky Blvd., 109028 Moscow, Russia}}
\email[]{altal (at) mi-ras (dot) ru }

\date{}

\begin{document}

\begin{abstract}
    In this note we generalise some of the work of Klarner on free semigroups of affine maps acting on the real line by using a classical approach from geometric group theory (the Ping--Pong lemma). We also investigate the boundaries within which Klarner's necessary condition for a semigroup to be related is applicable. 
\end{abstract}

\maketitle

\section{Introduction}

Several decades ago D.~Klarner studied a number of questions about integer affine maps that have connections to objects as diverse as self--orthogonal Latin squares and the famous Collatz conjecture \cite{Lagarias}. In particular, Klarner's results from~\cite{Klarner} that provide necessary and sufficient conditions for a collection of one--dimensional affine maps with integer coefficients to be generators of a free semigroup. A specific triple of such functions, namely $f_1(x)=2x$, $f_2(x)=3x+2$, and $f_3(x)=6x+3$, has raised the question if the density of the orbit of $x = 1$ under the iteration of $f_1$, $f_2$, and $f_3$ altogether is positive \cite[Problem 4]{Guy}. This problem still remains open.

This note shows that a few results from~\cite{Klarner} providing sufficient conditions for semigroups of integer affine maps to be free (and proved by using linear orders), can be deduced in a simple and  natural geometric way from the well--known ``Ping--Pong lemma'' in geometric group theory, when the extensions of these maps are considered to be acting on the real line. These results apply verbatim to the affine maps with real--valued coefficients, thus being independent of the number--theoretic context. In contrast, establishing that a given semigroup has non--trivial relations is a more delicate problem and is shown to be dependent on the number-theoretic nature of the coefficients. More information about semigroups arising from various actions of groups on the real line and circle can be found in the recent monograph~\cite{DNR}. 

Here and below, an affine map of $\mathbb{R}$ is a map $f: \mathbb{R} \longrightarrow \mathbb{R}$ defined by $f(x) = ax + b$, for some real $a, b \in \mathbb{R}$. Let $\mathrm{Aff}(\mathbb{R})$ be the group of all invertible affine maps relative to map composition (denoted by $\circ$ as usual). A semigroup $S$ generated by a set of maps $\{f_1,\ldots,f_n\}$ will be denoted by $S=\langle f_1, \ldots, f_r \rangle$. This semigroup is called \textit{free of rank} $r\geq 1$ with \textit{free basis} $f_1, \ldots, f_r$, if there are no two non-identical sequences $(i_1, \ldots, i_k)$ and $(j_1, \ldots, j_l)$ of indices in $\{1, \ldots, r\}$ such that $f_{i_1}\circ \ldots \circ f_{i_k} = f_{j_1} \circ \ldots \circ f_{j_l}$. Unless a different basis is explicitly indicated, we shall simply say that $S$ is a free semigroup assuming its free basis and rank.

A functional equality in $S = \langle f_1, \ldots, f_r \rangle$ of the form $f_{i_1}\circ \ldots \circ f_{i_k} = f_{j_1} \circ \ldots \circ f_{j_l}$, as above, will be referred to as a \textit{relation} in $S$. 

The Ping--Pong lemma can be formulated for semigroups of $\mathrm{Aff}(\mathbb{R})$ as follows (cf. \cite[VII.2]{delaHarpe} and also \cite[II.B]{delaHarpe} for a statement of the same flavour applied to free products of groups). It provides a simple but flexible condition for establishing a given semigroup be free (or a given group be free provided obvious modifications). 

\medskip
\noindent
\textbf{The Ping--Pong Lemma.} Let $S = \langle f_1, \ldots, f_r \rangle$, $r\geq 2$, be a semigroup of $\mathrm{Aff}(\mathbb{R})$. If there exists a collection of non-empty  mutually non-intersecting sets $I_1, \ldots, I_r$ such that $f_i(\cup_{j=1}^n I_j) \subset I_i$ for all $1\le i \le r$, then $S$ is a free semigroup of rank $r$ with free basis $f_1, \ldots, f_r$.
\medskip

\begin{proof}
Let $f_{i_1} \circ f_{i_2} \circ \ldots \circ f_{i_k} = f_{j_1} \circ f_{j_2} \circ \ldots \circ f_{j_l}$, with $k\geq l$, be a relation. Let $s \leq l$ be the maximal integer such that $(i_1, i_2, \ldots, i_s) = (j_1, j_2, \ldots, j_s)$. By multiplying the above relation on the left by $(f_{i_1} \circ f_{i_2} \circ \ldots \circ f_{i_s})^{-1}$ we cancel out the common prefix. Let $I = \cup_{j=1}^n I_j$. If $s < l$ then $f_{i_{s+1}}\circ \ldots \circ f_{i_k}(I) \subset I_{i_{s+1}}$ and $f_{j_{s+1}}\circ \ldots \circ f_{j_k}(I) \subset I_{j_{s+1}}$, while $I_{i_{s+1}} \cap I_{j_{s+1}} = \emptyset$, which is a contradiction. If $s = l$, then $f_{i_{s+1}}\circ \ldots \circ f_{i_k} = \mathrm{id}$, while $f_{i_{s+1}}\circ \ldots \circ f_{i_k}(I) \subset I_{i_{s+1}}\neq I$, which is also a contradiction. 
\end{proof}

The following result will follow almost immediately from the above version of the Ping--Pong lemma, and some elementary geometric considerations. 

\begin{theorem}\label{thm1}
Let $f_i(x) = a_i x + b_i$, for $i = 1, 2, \ldots, n$, be a finite subset of $\mathrm{Aff}(\mathbb{R})$ such that $a_i > 1$ for all $i = 1, 2, \ldots, n$. Let $s_i = \frac{b_i}{1-a_i}$. Then, up to a permutation of $f_i$'s, the semigroup $S = \langle f_1, \ldots, f_n \rangle$ is free whenever 
\begin{equation*}
    s_1 < s_2 < \ldots < s_n,
\end{equation*}
and
\begin{equation*}
    \frac{s_n - b_i}{a_i} \leq \frac{s_1 - b_{i+1}}{a_{i+1}},\quad \text{for all} \quad i = 1, 2, \ldots, n-1.
\end{equation*}
\end{theorem}

A quick comparison with \cite[Theorem 2.3]{Klarner} will convince the reader that the former and Theorem~\ref{thm1} are equivalent in the case of $a_i \geq 2$ and $b_i \geq 0$ being integers, for all $i = 1, 2, \ldots, n$ (with the order of these indices reversed). Otherwise, the number-theoretic nature of $f_i$'s is not essential in establishing that $S$ is a free semigroup: more details can be found in the proof of Theorem~\ref{thm1}. 

The following theorem generalises \cite[Theorem 2.2]{Klarner}. Its proof will follow from a basic picture illustrating the Ping--Pong Lemma (cf. Section 2).  

\begin{theorem}\label{thm-2-functions}
Let $f(x) = ax + b$, $g(x) = cx + d$ be two functions from $\mathrm{Aff}(\mathbb{R})$ such that $\sfrac{1}{a} + \sfrac{1}{c} \leq 1$. Then either $f$ and $g$ commute, or they generate a free semigroup. 
\end{theorem}

Again, here in Theorem~\ref{thm-2-functions} the number--theoretic nature of $f$ and $g$ is not important as compared to the proof provided in \cite{Klarner} to establish that the semigroup $S = \langle f, g \rangle$ is free. 

Speaking about the case when $S = \langle f_1, \ldots, f_n \rangle$ is not free, things become more subtle. Here the number--theoretic nature of $f_i$'s does manifest itself. First of all, we can easily show that \cite[Theorem 1.1]{Klarner} holds in a slightly more general context. 

\begin{theorem}\label{thm2}
Let $f_i(x) = a_i x + b_i$, for $i = 1, 2, \ldots, n$, be a finite subset of $\mathrm{Aff}(\mathbb{R})$ such that $a_i \in \mathbb{Z}$ are positive integers, and $b_i \in \mathbb{Q}$. Then the semigroup $S = \langle f_1, \ldots, f_n \rangle$ is not free (with free basis $f_1, \ldots, f_n$), whenever 
\begin{equation*}
    \frac{1}{a_1} + \ldots + \frac{1}{a_n} > 1.
\end{equation*}
\end{theorem}

However, obtaining a more general statement than Theorem~\ref{thm2} above appears difficult. The following examples are speaking to this matter. 

\begin{example}\normalfont
This example originates from a question of Erd\"{o}s that was solved by D.~J.~Crampin and A.~J.~W.~Hilton (see \cite[\S7]{Lagarias} for more details). If we take $f_1(x) = 2x + 1$, $f_2(x) = 3x + 1$, and $f_3(x) = 6x + 1$, then the respective semigroup $S = \langle f_1, f_2, f_3 \rangle$ is not free, as $f_1 \circ f_1 \circ f_2 = f_3 \circ f_1$. In this case, however, the inequality of Theorem~\ref{thm2} does not hold. 
\end{example}

\begin{example}\normalfont
Let $f_1(x) = 2x + 1$, $f_2(x) = 2x + \sqrt{2}$, and $f_3(x) = 2x + \sqrt{3}$. Then it is easy to check that all conditions of Theorem~\ref{thm2} are satisfied except for $b_i\in \mathbb Q$. However, the semigroup $S = \langle f_1, f_2, f_3 \rangle$ is free. Indeed, any relation in $S$ has the form \begin{equation*}
    f_{i_1} \circ f_{i_2} \circ \ldots \circ f_{i_L} = f_{j_1} \circ f_{j_2} \circ \ldots \circ f_{j_L},
\end{equation*}
for some $L\geq 2$, and implies a relation of the form 
\begin{equation*}
    \sum_{d\in I_1} 2^d + \sum_{d \in I_2} 2^d \sqrt{2} + \sum_{d \in I_3} 2^d \sqrt{3} = \sum_{d\in J_1} 2^d + \sum_{d \in J_2} 2^d \sqrt{2} + \sum_{d \in J_3} 2^d \sqrt{3},
\end{equation*}
for the constant terms of the functions, where $I_1 \cup I_2 \cup I_3 = J_1 \cup J_2 \cup J_3 = \{0, 1, 2, \ldots, L-1\}$, with $I_k$ indicating the instances of $f_k$ appearing on the left-hand side, and $J_k$ indicating the instances of $f_k$ appearing on the right-hand side of the relation above, for $k = 1, 2, 3$. The instances are numbered in the sequential order from left to right, starting from $0$ as the leftmost index. Thus, the amount of $f_k$ used and the order of composition can be completely recovered by ordering the indices in $\cup_{k=1,2,3} I_k$, respectively $\cup_{k=1,2,3} J_k$, in the decreasing order and then checking which index belongs to exactly which set $I_k$,  respectively $J_k$.   

From the fact that $1$, $\sqrt{2}$ and $\sqrt{3}$ are linearly independent over $\mathbb{Q}$, we obtain that $\sum_{d\in I_k} 2^d = \sum_{d\in J_k} 2^d$, for $k = 1, 2, 3$. The latter are just two binary presentations of the same number, which implies $I_k = J_k$, for $k = 1, 2, 3$. This is a contradiction since the index sequences have to be distinct.
\end{example}

\begin{example}\normalfont
Note that we cannot prove the statement of Theorem~\ref{thm2} for $a_i\in \mathbb Q$ either, as the maps $f_1(x) = x/2$ and $f_2(x) = (3x+1)/2$ generate a free semigroup, which follows from a general number--theoretical criterion \cite[Proposition~2]{Cas}. 
This semigroup is related to the Collatz conjecture \cite[\S 10]{Lagarias} and was also shown to be free in \cite[Theorem 2.2]{MisRod}, where the dynamical properties were investigated.
\end{example}

It is known that there is no algorithm which decides whether a given set of $3\times 3$ integer matrices represent a basis of a free semigroup \cite{KlaBiSa}, and this problem is undecidable even for upper--triangular matrices \cite{Cas}. This problem remains open for $2 \times 2$ integer matrices. The upper--triangular case, also considered in the work J.~Cassaigne, T.~Harju, and J.~Karhum\"aki \cite{Cas}, is easily seen to be equivalent to the case of groups generated by affine maps, via the correspondence \begin{equation*}
    ax + b \mapsto \left( \begin{array}{cc}
        a & b \\
        0 & 1
    \end{array} \right).
\end{equation*} 
In order to show that some of the semigroups considered  by in \cite{Cas} are free, they use some version of the Ping--Pong lemma (apparently without naming it) in conjunction with additional number--theoretic arguments.

Let us note that the group $\mathrm{Aff}(\mathbb{R})$ is \textit{metabelian} (solvable of derived degree $2$), and thus cannot contain any free subgroup of rank $r\geq 2$. However, by using advanced number--theoretical analysis one can prove that some groups generated by \textit{polynomial} maps are free, such as the group $\langle x \mapsto x+1, \,\, x \mapsto x^p \rangle$ with $p$ being an odd prime. This was first shown by S.~White \cite{White}. Later on, White's result was generalized in the works of S.~A.~Adeleke, A.~M.~W.~Glass, L.~Morley, and S.~D.~Cohen \cite{AGM, Cohen, CG}.

\section{Proofs}

\subsection*{Proof of Theorem~\ref{thm1}.} First, we pass to the inverses of the maps $f_i(x) = a_i x + b_i$, $i = 1, 2, \ldots, n$, which are given respectively by $g_i(x) = \frac{1}{a_i} x - \frac{b_i}{a_i}$, $i = 1, 2, \ldots, n$. It is clear that $S$ is a free semigroup if and only if $\bar{S} = \langle g_1, \ldots, g_n \rangle$ is free.

\begin{figure}[h]
    \centering
    \includegraphics[scale=0.20]{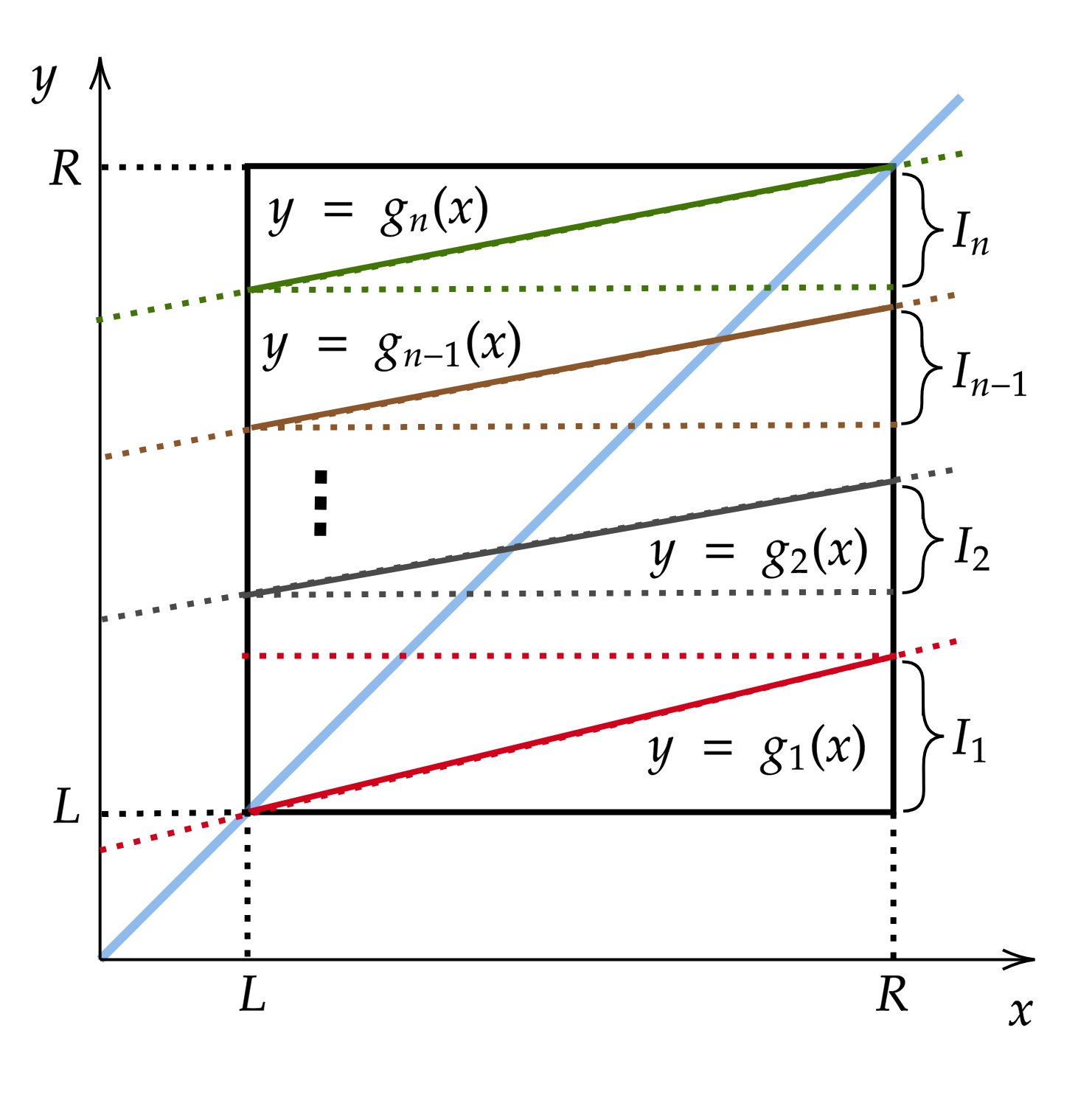}
    \caption{The graphs of the functions $y = g_1(x)$, $\ldots$, $y = g_n(x)$ in the square $B = [L,R]\times [L,R]$. The intervals $I_1, \ldots, I_n \subset [L, R]$ are marked on the right}
    \label{fig:intervals}
\end{figure}

The next geometric considerations are shown in Figure \ref{fig:intervals}. Consider the graphs of the functions $g_1(x), \ldots, g_n(x)$ and the graph of the identity (or ``diagonal'') function $D(x)=x$ on the plane $\mathbb R^2$ (see Figure~\ref{fig:intervals}). It is easy to compute that the graphs of $g_i(x)$ and $D(x)$ intersect in the point $u_i=(s_i,s_i)$ for $s_i=\frac{b_i}{1-a_i}$. Since we have $s_1<s_2<\ldots<s_n$, the intersection points $u_i$ are located on the diagonal according to their indices. Let us put $L = s_1, R = s_n$ and let $B$ denote the square $B = [L, R]\times [L, R]$.

The following conditions are just Theorem's premises being reformulated:

\begin{equation}\label{LR}
    g_1(L) = L, \quad g_n(R) = R,\quad L < R,
\end{equation}
\begin{equation}\label{intervals1}
    g_i(R) \leq g_{i+1}(L), \text{ for all } i = 1, \ldots, n-1.
\end{equation}
The second condition means geometrically that the intersection of $y = g_i(x)$ with the right---hand side of the square $B$ is below the intersection of $y = g_{i+1}$ with the left---hand side of $B$. Since both functions are increasing, it also means that their restrictions to the interval $I=(L,R)$ have non--intersecting images $I_i=g_i(I)$ and $I_{i+1}=g_{i+1}(I)$, with $I_{i}$ being located strictly below $I_{i+1}$. Hence, all images $I_1, I_2, \ldots, I_n$ are located on the right---hand side of the square $B$, from bottom up to the top, and do not intersect. Then the Ping--Pong lemma applies to $\bar{S}$ with the collection of intervals $\{I_i\}^n_{i=1}$. \qed

\subsection*{Proof of Theorem~\ref{thm-2-functions}.} The proof follows from Figure~\ref{fig:intervals}, where only the lowest and the upmost functions are present. In the notation of Theorem \ref{thm1}, let $g_1 = f^{-1}$ and $g_2 = g^{-1}$. The intervals $I_f = I_1$ and $I_g = I_2$ touch the left and right end of the interval $(L, R)$ correspondingly, and since $\sfrac{1}{a} + \sfrac{1}{c} \leq 1$, we have $|I_1|+|I_2|\leq R-L$, hence $I_1$ and $I_2$ do not intersect. The Ping--Pong lemma can be applied to show that the semigroup $\langle f, g \rangle$ is free, unless the intervals $I_1$ and $I_2$ degenerate in the case $L = R$, which is equivalent to $f$ and $g$ commuting, since this condition can be expressed as $f(g(0)) = g(f(0))$. \qed

\subsection*{Proof of Theorem~\ref{thm2}.} The scheme of the proof is the same as in \cite[Theorem 1.1]{Klarner} for affine functions with integer coefficients. First, the general case is reduced to the case of $a_1 = \ldots = a_n = m$, with $n > m$, and then the theorem is proved under the latter assumption. (The proof of the latter part in \cite{Klarner} is omitted and attributed to R.~Rivest, but the suggested reference leads to the unpublished manuscript \cite{Klarner-man} instead). 

\medskip

The reduction to the case $a_1 = \ldots = a_n = m$, with $n > m$ in \cite[Theorem 1.1]{Klarner} works also for rational coefficients $b_i$, and we include this argument for completeness. Let $S_L \subset S$ be the set of all possible (distinct) length $L$ words over the alphabet $\Sigma = \{f_i\}^n_{i=1}$ (i.e. formal compositions of functions from the generating set). 

We can also consider the elements of $S_L \subset S$ as $L$--fold compositions of affine functions. For $f \in S$ of the form $f(x) = a\,x + b$ we have $f'(0) = a$. Thus,
\begin{equation*}
    \sum_{f \in S_L} \frac{1}{f'(0)} = \left( \frac{1}{a_1} + \frac{1}{a_2} + \ldots + \frac{1}{a_n} \right)^L = \mu^L.
\end{equation*}
The latter can be also written as a multinomial sum
\begin{equation*}
    \left( \frac{1}{a_1} + \frac{1}{a_2} + \ldots + \frac{1}{a_n} \right)^L = \sum_{l_1 + \ldots + l_n = L} \binom{L}{l_1, l_2, \dots, l_n}\, a_1^{-l_1}\, a_2^{-l_2} \ldots a_n^{-l_n}.
\end{equation*}
Let us consider a sufficiently large $L$ for which $\mu^L > L^n$. This is indeed possible due to the fact that $\mu>1$. For each $L$ the number of non-negative tuples $(l_1, l_2, \dots, l_n)$ with $l_1 + l_2 + \dots + l_n = L$ does not exceed $L^n$, which implies that there exists a tuple $(l_1, \ldots, l_n)$ such that 
\begin{equation*}
    \binom{L}{l_1, l_2, \dots, l_n}\, a_1^{-l_1}\, a_2^{-l_2} \ldots a_n^{-l_n} \geq \frac{\mu^L}{L^n} > 1.
\end{equation*}
\smallskip

Let $E$ be the set of elements $f \in S_L$ with $f'(0) = m = a_1^{l_1} a_2^{l_2} \ldots a_n^{l_n}$. Then $|E|\geq\binom{L}{l_1, l_2, \ldots, l_n} > m$. If two elements in $E$ coincide as affine functions, then $S$ has a relation. 
Otherwise, we can replace the initial semigroup $S$ by $\langle E \rangle$ and assume that all $f_i = a_i x + b_i$ satisfy $a_i = m$, $i = 1, \ldots, n$, with $n > m$.    
\medskip

Now we generalize the argument of R.~Rivest to the case when $a_i\in \mathbb Z$ and $b_i \in \mathbb Q$. For each function $f_i(x) = a_i x + b_i$ in $\Sigma = \{f_i\}^n_{i=1}$, let $b_i = \frac{p_i}{q_i}$, $i = 1, 2, \ldots, n$. Note that we have $a(L) = |S_L| = n^L$. On the other hand, each function $f \in S_L$ has the form 
\begin{equation*}
    f(x) = m^L x + \sum^{L-1}_{d = 0} m^d\, \frac{p_{i_d}}{q_{i_d}}.
\end{equation*}
Since each constant term can be written as some fraction $A+B/C$ subject to the inequalities $|A|\leq (L m^L \max_{i=1,\ldots,n} \lceil b_i \rceil)$ and $0<|B|<C\leq \mathrm{lcm}\{ q_i \}^n_{i=1}$, we obtain that
the number of distinct constant terms is bounded above by the quantity 
\begin{equation*}
    b(L) = L \cdot m^L \cdot \max_{i=1,\ldots,n} \lceil b_i \rceil \cdot \mathrm{lcm}\{ q_i \}^n_{i=1}.
\end{equation*}

For $L$ large enough, $n > m$ implies that $a(L) > b(L)$, and thus we have two distinct words over $\Sigma$ that coincide as functions in $\mathrm{Aff}(\mathbb{R})$. \qed 
\section*{Acknowledgements}

The authors would like to thank Jeffrey~C.~Lagarias and the anonymous referee for their comments and literature suggestions. A.K. was partially supported by SNSF (project no.~PP00P2-170560) and Russian Federation Government (grant no. 075-15-2019-1926). A.T. was partially supported by RFBR and SC RA (project no.~20-51-05010). Also, A.T. would like to thank the University of Neuch\^atel for its hospitality in May 2021.

\end{document}